\newif\ifSV
\newif\ifLNCS
\SVfalse
\LNCSfalse
\ifSV
\documentclass{svjour3}
\else
\ifLNCS
\documentclass{llncs}
\else
\documentclass{article}
\fi
\fi
\usepackage{a4}
\usepackage[all]{xy}
\usepackage{qsymbols}
\usepackage{marvosym}
\usepackage{fancybox}
\usepackage{url}
\usepackage{color}
\usepackage{haskell}

\newcommand{\nat}{\ensuremath{\mathbb{N}}}
\newcommand{\real}{\ensuremath{\mathbb{R}}}

\newcommand{\fl}[1]{\ar@/^/[#1]^r \ar@/^/[d]^{d}}
\newcommand{\flr}[1]{\ar@[blue]@2@/^/[#1]^{\color{blue} r} \ar@/^/[d]^{d}}
\newcommand{\fld}[1]{\ar@/^/[#1]^r \ar@[blue]@2@/^/[d]^{\color{blue} d}}

\newcommand{\Ealg}{eAlg}
\newcommand{\Galg}{gAlg}

\title{An exercise on streams: convergence acceleration} %
\ifSV 
\author{Pierre Lescanne\\
  University of Lyon, \\
  \'Ecole normale sup\'erieure de Lyon, \\
  LIP (UMR 5668 CNRS ENS Lyon UCBL INRIA)\\
  46 all\'ee d'Italie, 69364 Lyon, France} 
\else 
\ifLNCS
\author{Pierre Lescanne}
\institute{University of Lyon, \\
  \'Ecole normale sup\'erieure de Lyon, \\
  LIP (UMR 5668 CNRS ENS Lyon UCBL INRIA)\\
  46 all\'ee d'Italie, 69364 Lyon, France} 
\else
\author{Pierre Lescanne\\
  University of Lyon, \\
  \'Ecole normale sup\'erieure de Lyon, \\
  LIP (UMR 5668 CNRS ENS Lyon UCBL INRIA)\\
  46 all\'ee d'Italie, 69364 Lyon, France} 
\fi
\fi
\date{}

\begin{document}
\maketitle

\begin{abstract}
\ifSV\else \medskip  \hrule\fi

\medskip

This paper presents convergence acceleration, a method for computing efficiently the
limit of numerical sequences as a typical application of streams and higher-order
functions.

\ifSV\end{abstract}


\else \medskip

\end{abstract}
\hrule\fi

\bigskip 

\section{Introduction}
\label{sec:intro}

Assume that we want to compute numerically the limit of a sequence that converges
slowly.  If we use the sequence itself, we will get significant figures of the limit
after a long time. Methods called \emph{convergence acceleration} have been designed
to shorten the time after which we get reasonable amount of significant digits.  In
other words, \emph{convergence acceleration} is a set of methods for numerically
computing the limit of a sequence of numbers.  Those methods are based on sequence
transformations and are a nice domain of application of
streams~\cite{DBLP:journals/tcs/Rutten03,DBLP:conf/cmcs/BonsangueRW12}, with
beautiful higher order functions.  This allows us to present elegantly rather complex
methods and to code them in \textsf{Haskell}~\cite{hutton07:_progr_haskell} replacing
long, obscure and special purpose \textsf{Fortran} programs by short, generic and
arbitrary precision \textsf{Haskell} programs.

In this paper we show, how given a sequence $(s_n)_{n`:\nat}$, we can evaluate
efficiently $\displaystyle \lim_{n"->"\infty} s_{n}$. For that we use \emph{Levin
  transforms}. There are three kinds of such transforms, which are the result of
three sequence transformations labeled traditionally by $t$, $u$ and
$v$~(\cite{osada93:_accel_method}, p.58).

\section{Presentation of the method}
\label{sec:pres}

\nocite{opac_b1101304}
\nocite{laurie04:_conver}

In what follows we speak indistinctly of ``sequences'' or of ``streams'' We use
\textsf{Haskell} and we work with arbitrary precision reals based on the
implementation of David Lester called $CReal$.  In \textsf{Haskell} the type
of streams over a type $A$ is written~$[A]$.

For the numerical aspect, we follows Naoki Osada~\cite{osada93:_accel_method} (for a
historical account see~\cite{Brezinski00convergenceacceleration}) and we
show that the stream notation of \textsf{Haskell} makes the presentation
simpler. With no surprise, the size of the \textsf{Haskell} code is the same if not
shorter than the mathematical description of the
reference~\cite{osada93:_accel_method}. Moreover it provides efficient programs.

Levin transformations are somewhat generic in the sense that they are based on elementary
transformations. Specialists of convergence acceleration propose three such elementary
transformations.  Let $s$ be a sequence on $CReal$, i.e., \<s :: [CReal]\>. We define
first a basic sequence transformation on which we will found our elementary transformations:
\begin{haskell}
dELTA :: [CReal] "->" [CReal]\\
dELTA s  =  zipWith (-) (tail s) s
\end{haskell}
which means that $dELTA(s)_n = s_{n+1} -s_n$.  From this basic sequence
transformation we define the three elementary other sequence transformations as
follows.  A unique function depending on a parameter which is either \<T\> or
\<U\> or \<V\> is given.  For a computer scientist those names \<T\>,
  \<U\>, \<V\> are meaningless but this terminology (in lower case i.e., $u$,
  $t$, $v$ is traditionally used by mathematicians and we stick to it. It
corresponds to the traditional notations of numerical analysis for those sequence
transformations.
\begin{haskell}
data Kind = T | U | V\\~\\
delta :: Kind "->" [CReal] "->" [CReal]\\
delta T s  =  dELTA s\\
delta U s  =  zipWith (*) (dELTA s) [1..]\\
delta V s  =   zipWith\ (/)\ \hsalign{(zipWith\ (*)\ (tail\$dELTA s)\ (dELTA s))\\
                             (dELTA (dELTA s))}
\end{haskell}

In numerical analysis, people speak about \emph{E-algorithm}. This is a family of functions
$\Ealg_{n,k}$ which  are also parametrized by a character  either \<T\>
or \<U\> or \<V\>. It tells which of the basic sequence transformations is
chosen. $\Ealg_{n,k}$  uses a family of auxiliary functions which we call
$\Galg_{n,k}$  for symmetry and regularity.
Here is the \textsf{Haskell} code for these functions:
\begin{haskell}
eAlg :: Kind "->" Int "->" [CReal] "->" [CReal]\\
eAlg c 0 s = s\\
eAlg c k s = \hslet{ 
                   a = (eAlg c (k-1) s)\\
                   b = (gAlg c (k-1) k s)
                   }{%
                   zipWith (-) a (zipWith (*) b (zipWith (/) (dELTA a) (dELTA b)))
                   }
\end{haskell}

\vspace{5pt}

\begin{haskell}
gAlg  :: Kind "->" Int "->" Int "->" [CReal] "->" [CReal]\\
gAlg c 0 j s = \hslet{
                 nTojMinus1 j = zipWith (**) [1..]  (repeat (fromIntegral (j-1)))
                 }{%
                 zipWith (/) (nTojMinus1 j) (delta c s)
                 }\\
gAlg c k j s = \hslet{
                 a = gAlg c (k-1) j s\\
                 b = gAlg c (k-1) k s
                 }{%
                 zipWith (-) a (zipWith (*) b (zipWith (/) (dELTA a) (dELTA b)))
                 }
\end{haskell}
Here is the formula as it is given in~\cite{osada93:_accel_method}.  $R_n$ is the
generic value of \<(delta c s)_n\>. \<gAlg\> is written $g$. 
$E_k^{(n)}$ is the $n^{\rm th}$ element of the sequence
\<eAlg c k s\>, the same for $g_{k,j}^{(n)}$. $`D$ is the notation for what we write
  \<dELTA\>.
  \begin{displaymath}
    \begin{array}{lcl@{,\qquad}l@{;~}l}
E_0^{(n)} &=& s_n,\ \ g_{0,j}^{(n)} = n^{1-j} R_n & n =1,2,... & j = 1,2,...,\\[5pt]
E_k^{(n)} &=& E_{k-1}^{(n)} - g_{k-1,k}^{(n)}
\frac{\displaystyle`DE_{k-1}^{(n)}}{\displaystyle`Dg_{k-1,k}^{(n)}}&n=1,2,...&k=1,2,...,\\[8pt]     
g_{k,j}^{(n)} &=& g_{k-1,j}^{(n)} - g_{k-1,k}^{(n)}
\frac{\displaystyle`Dg_{k-1,j}^{(n)}}{\displaystyle`Dg_{k-1,k}^{(n)}}&n=1,2,...&k=1,2,..., j>k
    \end{array}
  \end{displaymath}
\section{Levin's formulas}
\label{sec:levin}
There is another formula for the $E-algorithm$:
\begin{displaymath}
  T_k^{(n)} = \frac{\displaystyle\sum_{j=0}^k (-1)^j {k \choose j}{n+j \choose
      n+k}^{k-1} \frac{s_{n+j}}{R_{n+j}}}
{\displaystyle\sum_{j=0}^k (-1)^j {k \choose j}{n+j \choose n+k}^{k-1} \frac{1}{R_{n+j}}}
\end{displaymath}
$T_k^{(n)}$ is not easily amenable to a Haskell program.\footnote{In particular due to many
  divisions by $0$ and to the complexity of the formula.}
We give only the functions for $k=0, 1, 2$, which we call \<levin\ c\ 0\>, \<levin\ c\ 1\> and
\<levin\ c\ 2\>. For $k=0$, $T^{(n)}_0 = s_n$ and \<levin k\ 1\> and \<levin k\ 2\> are the result of small
calculations. 
\begin{haskell}
levin  :: Kind "->"  Int "->" [CReal] "->" [CReal]\\
levin c 0 s  =  s\\
levin c 1 s  =  zipWith (-) s (zipWith (/) \hsalign{(zipWith (*) (dELTA s) (delta c s)) \\
                                                     (dELTA (delta c s)))}
\end{haskell}
\<levin U 1 s\> is called \emph{Aitken's delta-squared process}.
One notices that the numerator and the denominator of the above formula differ
slightly. Indeed $s_{n+j}$ in the numerator is just replaced by $1$ in the
denominator.
\begin{haskell}
\hspace*{-30pt}formulaForLevinTwo :: Kind "->" [CReal] "->" [CReal] "->" [CReal]\\
\hspace*{-30pt}formulaForLevinTwo c s' s = \\
\qquad   \hsalign{zipWith (+) \hsalign{(zipWith (-) \hsalign{(foldl (zipWith (*)) [2..] [tail\$tail s',tail\$delta c s,delta c s])\\
                                                                       (foldl (zipWith (*)) [2,4..] [tail s',tail\$tail\$delta c s,delta c s]))}\\
                                         (foldl (zipWith (*)) [0..] [s', tail\$tail\$delta c s,tail\$delta c s])}}
\end{haskell}
\begin{haskell}
\hspace*{-30pt}levin c 2 s = zipWith (/) (formulaForLevinTwo c s s) (forLevinTwo c [1,1..] s)\\
\end{haskell}
Brezinski~\cite{brezinski80} proves that the sequences $E_k^{(n)}$ and $T_k^{(n)}$ are
the same, in other words:
\begin{haskell}
eAlg =_{`h} levin
\end{haskell}

\section{Classic and non classic examples}
\label{sec:examples}

\rightline{\parbox{8cm}{\begin{it} More than any other branch of numerical analysis,
      convergence acceleration is an experimental science. The researcher applies the
      algorithm and looks at the results to assess their worth. 
\end{it}
}}
\medskip
\rightline{Dirk Laurie~\cite{laurie04:_conver}}

\bigskip

Given a sequence we use the convergence acceleration to find the main coefficient of
the asymptotic equivalent.  More precisely given an integer sequence $s_n$, we want to find
a number $a >1$ such that $s_n\sim a^n f(n)$ where $f(n)$ is subexponential.  In other
words, for all $b`:\real$, $b>1$, $\displaystyle \lim_{n"->"\infty} \frac{f(n)}{b^n}
= 0$. Actually we compute the limit of the sequence $s_{n+1}/s_n$.  For that we create the
function:
\begin{haskell}
  expCoeffAC :: ([CReal] "->" [CReal]) "->" [Integer] "->" Int "->" CReal\\
  expCoeffAC transform sequence n =   last (transform  (zipWith (/) (tail u) u))
             \hswhere{u = map fromIntegral (take n sequence)}\\
\end{haskell}
Thus \<expCoeffAC (levin U 2) s 300\> gives the approximation of the
coefficient one can get after $300$ iterations using the sequence transformation
\<levin U 2\>.

\subsection{Catalan numbers}
\label{sec:catalan}

A good example to start is Catalan numbers:
\begin{haskell}
  catalan = 1 : [\hsalign{\textbf{let} cati = take i catalan\\
  \textbf{in} sum (zipWith (*) cati (reverse cati)) | i "<-" [1..]]}
\end{haskell}
We know that 
\begin{displaymath}
  catalan !! n \sim \frac{4^n}{\sqrt{\pi n^3}}
\end{displaymath}
Actually we get \< expCoeffAC (levin U 2) catalan 800 = 4.0000000237\> (with $8$
exacts digits) and
\begin{eqnarray*}
expCoeffAC~(eAlg~T~2)~catalan~800 &\approx&
3.9849561088\\
  expCoeffAC~(eAlg~U~2)~catalan~800 &\approx&
3.9773868157\\
 expCoeffAC~(eAlg~V~2)~catalan~800 &\approx&
3.9773869346
\end{eqnarray*}

\subsection{Counting plain lambda terms}
\label{sec:count-lmbd}

Now we want to use this technique to address a conjecture\footnote{This problem is
  the origin of the interest of the author for this
  question.}~\cite{GrygielLescanne-Binary}, on the asymptotic evaluation of the
exponential coefficient of the numbers of typable terms of size $n$ when $n$ goes
to~$\infty$.  First let us give the recursive definition of the numbers $S_{\infty}$
of plain lambda terms of size $n$. This sequence appears on the \emph{On-line
  Encyclopedia of Integer Sequences} with the entry number \textbf{A114851}.  We
assume that abstractions and applications have size $2$ and variables have size $1+k$
where $k$ is the depth of the variable w.r.t. its binder.
\begin{eqnarray*}
  S_{\infty,0} &=& S_{\infty,1} ~=~ 0,\\
  S_{\infty,n+2} &=& 1 + S_{\infty,n} + \sum_{k=0}^n S_{\infty,k} S_{\infty,n-k}.
\end{eqnarray*}
It has been proved in~\cite{GrygielLescanne-Binary}  that
\[ S_{\infty,n} \sim A^n \cdot \frac{C}{n^{3/2}},\] where $A \doteq
  1.963447954$ and $C \doteq 1.021874073$.  After $300$ iterations and using \<levin
  U 2\> we found \[\mathbf{1.96344}89522735283291619147713569993355616.\]
giving six exact digits. 

\subsection{Counting typable lambda terms}
\label{sec:count-typ-lmbd}

The question is now to find the exponential coefficients for the numbers of typable
terms.  We have no formula for computing those numbers. The only fact we know is the
following table of the numbers $T_{\infty,n}$ till $42$ which has been obtained after heavy computations
(more than 5 days for the $42^{\rm nd}$).  The method consists in generating all the lambda terms
of a given size and sieving those that are typable to count them.
Therefore the best method to guess the exponential coefficient is by acceleration of
convergence.  After $43$ iterations we found $1.8375065809...$. Knowing that with the
same number $43$ we get $1.8925174623...$ for $S_{\infty}$, this is not enough to
conclude. But this allows us to speculate that the exponential coefficient for the
asymptotic evaluation of $T_{\infty,n}$ could be $1.963447954$ like for
$S_{\infty,n}$'s.

\section{Application to divergent series}
\label{sec:Euler}

In his famous paper~\cite{euler60:_de} Euler provides a sum to divergent series.
See~\cite{ramis93:_séries,martinet91:_elemen_accel_multis} and for a light introduction,
the reader who understands French is advised to watch
the video~\cite{ramis09:_leonh_euler} which completes another
video~\cite{numberphile14:_astoun} in English.  

Among the methods Euler and his followers propose to give a meaning to sum of divergent series there
is convergence acceleration.  We applied naturally our implementation to some
divergent series.  

Lets us define the function that associated to a sequence it series.
\begin{haskell}
seq2series :: [Int] "->" [CReal]\\
seq2series s = \hsalign{\textbf{let} series (x:s') = x : map ((+) x) (series s')\\
               \textbf{in} map fromIntegral (series s)}
\end{haskell}

\subsection{Grandi series }
\label{sec:1+1+1}
Grandi series is also called by Euler, Leibniz series. This is the series
\begin{displaymath}
  \sum_{i=0}^\infty (-1)^{i}
\end{displaymath}
which is sometime written $1-1+1-1+\dots$.  In \textsf{Haskell} it is:
\begin{haskell}
grandi = \textbf{let} gr = 1:(-1):gr \textbf{in} seq2series gr
\end{haskell}
We get $1/2$ after $3$ iterations using \<eAlg T 2\> and this does not change when we
increase the number of iterations.

\subsection{Leibniz series 1 - 2 + 3 - 4 + 5 + ...}
This series attributed to Leibniz is also studied by Euler:
\begin{displaymath}
  \sum_{i=0}^\infty (-1)^{i+1} i
\end{displaymath} 
In \textsf{Haskell}
\begin{haskell}
sumNatAlt = \textbf{let} lbn i = -i : (i+1) : lbn (i+2) \textbf{in} seq2series (lbn 0)
\end{haskell}
By using \<eAlg U 4\> we get $0.25$ after $6$ iterations in accordance with Euler's result~$\frac{1}{4}$.

\section{Conclusion}

We have shown how streams can be applied to numerical analysis, namely to convergence
acceleration.  It makes no doubt that they can also be applied to other fields in
numerical analysis or elsewhere. For instance, one may imagine applications of
acceleration of convergence to the computation of limits of non numerical sequences.

We did not use the differential
equation approach~\cite{DBLP:journals/tcs/Rutten03,DBLP:conf/calco/Winter13}, but presenting
acceleration of convergence in this framework should be worthwhile. 


\end{document}